\documentclass[twoside,12pt]{amsart} 
\usepackage{amsfonts,amssymb,latexsym,amsthm, epsfig,amsmath,amscd}
\newcommand{\mf}{\mathfrak} 

\newcommand{\inj}{\hookrightarrow} 
\newcommand{\lra}{\longrightarrow}

\newcommand{\ol}{\overline}

\newcommand{\ra}{\rightarrow}
\newcommand{\Ra}{\Rightarrow}

\newcommand{\eps}{\varepsilon}
\newcommand{\stk}{\stackrel}

\newcommand{\mbb}{\mathbb}
\newcommand{\tn}{\textnormal}

\newtheorem{de}{Definition}[section]
\newtheorem{re}[de]{Remark}
\newtheorem{pr}[de]{Proposition} 
\newtheorem{tr}[de]{Theorem}
\newtheorem{lm}[de]{Lemma} 

\newtheorem{nt}[de]{Notation} 

\newtheorem{co}[de]{Corollary}

\newcommand{\LMD}{\Lambda}
\newcommand{\lmd}{\lambda}
\newcommand{\lb}{\linebreak}
\def\vp{\rm \vspace{0.2cm}}
\def\M{\rm Max}
\def\m{\mf{m}}

\def\hb{\hfill$\Box$}

\def\Um{\rm Um}
\def\GL{\rm GL}
\def\GH{\rm GH}
\def\GQ{\rm GQ}
\def\SL{\rm SL}
\def\EQ{\rm EQ}

\def\EH{\rm EH}
\def\SH{\rm SH}
\def\SQ{\rm SQ}

\def\E{\rm E}
\def\EH{\rm EH}

\def\G{\rm G}
\def\GQ{\rm GQ}
\def\EM{\rm M}
\def\Sp{\rm Sp}

\def\k{\rm K_1}
\def\K{\rm K}
\def\KH{\rm KH}

\def\I{\rm I}

\def\O{\rm O}

\def\C{\rm C}
\def\es{\rm S}

\begin{document}
\title{Local-global principle for General Quadratic and General Hermitian groups
  and the nilpotence of ${\KH}_1$}
\author{Rabeya Basu}
\address{Indian Institute of Science Education and Research - Pune,
Maharashtra 411008, India}
\email{rabeya.basu@gmail.com, rbasu@iiserpune.ac.in}

\date{}
\maketitle
%\begin{center} 
%{\it A tribute to the late Professor Amit Roy}
%\end{center} \vp 

\begin{center} 
{\bf Abstract}
\end{center} 

\noindent {\small  In this article we establish an analog of the Quillen---Suslin's 
local-global principle for the elementary subgroup of the general quadratic group and 
the general Hermitian group. We show that unstable ${\k}$-groups 
of general Hermitian groups over module finite rings are nilpotent-by-abelian. This generalizes earlier results of 
A. Bak, R. Hazrat, and N. Vavilov.} \vp 

\section{\large Introduction}

The vigorous study of general linear groups and more generally algebraic ${\K}$-theory was 
stimulated in mid-sixties by the desire to solve Serre's problem 
on projective modules ({\it cf.} Faisceaux Alg\'ebriques Coherent, 1955). This prominent problem in commutative algebra asks whether 
{\it finitely generated projective modules over a polynomial ring over a field are free}. 
The beautiful book 
{\it Serre's problem on projective modules} by T.Y. Lam  gives a comprehensive account 
of the mathematics surrounding Serre's problem and its solution. Later we see analogs of Serre's 
Problem for modules with forms and for other classical groups in the work of H. Bass, A. Suslin, L.N. Vaserstein,
V.I. Kopeiko, R. Parimala and others in \cite{Bass}, \cite{SUSV}, \cite{KOP}, \cite{SUSK}, \cite{P1}, \cite{P2}. 
In this current paper, we are interested in the context of modules with forms in certain problems related to Serre's Problem, {\it viz.} 
normality  of the elementary subgroup of the full automorphism group, 
Suslin's local-global principle for classical-like groups, stabilization for ${\k}$-functors of classical-like groups, and the 
structure of unstable ${\k}$-groups of classical-like groups. 

Difficulties one has in handling the 
quadratic version of Serre's Problem in characteristic 2 were first 
noted by Bass in \cite{Bass}. In fact, in many cases it was difficult 
to handle classical groups over fields of characteristic 2, rather 
than classical groups over fields of char $\ne$ 2. 
(For details see \cite{HV1}). In 1969, A. Bak resolved this problem by 
introducing {\it form rings} and {\it form parameter}. He 
introduced the general quadratic group or Bak's unitary group, 
which covers many different types of classical-like groups.  
We also see some results in this direction in the work of 
Klein, Mikhalev, Vaserstein {\it et al.} in 
\cite{K1}, \cite{K2}, \cite{V2}. The concept of {\it form parameter} also 
appears in the work of K. McCrimmon, and plays an important role in his 
classification theory of Jordan algebras ({\it cf.}~\cite{M}), 
for details see (\cite{HO}, footnote pg. 190.) and  \cite{J}. 
In his seminal work ``${\K}$-theory of forms", Bak has established 
analog of many problems related to Serre's problem in a very explicit and 
rigorous manner. 
But, Bak's definition of the general quadratic group does not 
include many other types classical-like groups, {\it viz.} 
odd dimensional orthogonal groups, exceptional groups  ${\E}_6$, ${\E}_7$, ${\E}_8$ {\it etc.} 
In 2000, G. Tang, in his Ph.D thesis, established  analog of many results
for the general Hermitian groups. 
Very recently, in 2005, Victor Petrov using Bak's concept of 
{\it doubly parametrized form parameter}
has resolved this problem by introducing {\it odd unitary groups}, 
which also includes Bak's unitary and general Hermitian groups; {\it cf.~\cite{P}}.  
Also, he has established many analogous results  for his group.

In 1976, D. Quillen came up with a localization method 
which was one of the main ingredients for the proof of 
Serre's problem (now widely known as Quillen---Suslin Theorem).
Shortly after the original proof Suslin introduced the 
following matrix theoretic version of Quillen's local-global 
principle. 

{\bf Suslin's Local-Global Principle:} 
{\it Let $R$ be a commutative ring with identity, $X$ a variable and 
$\alpha(X)\in {\GL}(n, R[X])$ with $\alpha(0)={\I}_n$, $n\ge 3$. If 
$\alpha_{\mf{m}}(X)\in {\E}(n, R_{\mf{m}}[X])$ for every maximal ideal 
$\mf{m}\in {\M}(R)$, then $\alpha(X)\in {\E}(n, R[X])$. } \vp

Soon after he gave the ${\k}$-analog of Serre's problem, 
which says, \vp

{\it for a polynomial ring in $r$ variables over a field $K$ elementary subgroup 
of ${\GL}_n(R)$ coincides with the special linear group. i.e.}
$${\E}_n(K[X_1,\ldots,X_r])={\SL}_n(K[X_1,\dots,X_r]). $$
In connection of this theorem 
he proved the normality of the 
elementary subgroup ${\E}(n, A)$ in  the general linear 
group ${\GL}(n, A)$, over a module finite ring $A$, for $n\ge 3$; ({\it cf.}~\cite{Tu}). 
Later analogous results for the symplectic and orthogonal groups 
were proven by Suslin and Kopeiko in \cite{SUS} and 
\cite{SUSK} and by Fu An Li in \cite{Fu}, and for arbitrary Chevalley 
groups by Abe ({\it cf.}~\cite{A}) in the local case, and 
by Taddei  ({\it cf.}~\cite{Ta}) in general. Later we see a simpler 
and more general treatment in  works of Ambily, Bak, Hazrat, Petrov, Rao, Stavrova, Stepanov, Suzuki, Vavilov, and others.

We see generalization of the above local-global principle for the symplectic group in \cite{KOP}, and for 
the orthogonal group in \cite{SUSK}. 
The normality of the general quadratic groups is known from the 
work of A. Bak and N. Vavilov, {\it cf.\cite{BV}}. In \cite{T}, 
G. Tang has proved the normality property for the general 
Hermitian groups. 
In \cite{BRK},  we have shown that the question of normality of the
elementary subgroup of the 
general linear group, symplectic and orthogonal groups, is equivalent 
to the above local-global principle, where the base ring is associative with identity and 
finite over its center. In that article above three classical groups were treated uniformly. 
Motivated by the work of 
A. Bak, R.G. Swan, L.N. Vaserstein and others,  in \cite{BBR}, the author with 
A. Bak and R.A. Rao has  established an analog
of Suslin's local-global principle for the transvection subgroup of 
the automorphism group of projective, 
symplectic and orthogonal modules of global rank at least 1 and 
local rank at least 3, under the assumption that the projective module has 
constant local rank and that the symplectic and orthogonal modules are locally 
an orthogonal sum of a constant number of hyperbolic planes.
In this article we have proved the equivalence of the local global principle 
with the normality property. Since normality holds in the above cases, 
this establishes that the local global principle also holds. 
In fact following Suslin-Vaserstein's method we establish an analogous local-global principle 
for the general quadratic and general Hermitian groups. 

We treat these two groups uniformly and give explicit proofs of those results.
We have overcome many technical difficulties which come in the 
Hermitian case due to the elements $a_1,\ldots,a_r$ 
(with respect to these elements we define the Hermitian groups). 
We assume $a_1=0$. The rigorous study of the general Hermitian groups can be found in \cite{T}. 
In \cite{BV}, we get an excellent survey on this area 
in a joint work of A. Bak and N. Vavilov. We refer \cite{HVZ} for an alternative approach to
localization, \cite{HSVZ}  for a general overview, and \cite{CR} for relative cases. 
Also, for commutative rings with identity
the Quillen---Suslin's local-global principle is in the work of V. Petrov and A. Stavrova ({\it cf.}~\cite{PS}), which
covers, in particular, classical groups of Witt index $\ge 2$ or $\ge 3$, depending on the type.

In \cite{BRK}, it has been shown that the normality criterion 
of the elementary subgroup of the general linear group 
is equivalent to the above local-global principle. In this paper 
we establish the analogous local-global principle for the general 
quadratic and Hermitian group, and prove an 
equivalence.  More precisely, we prove ($\S 6$, Theorem \ref{LG},  and $\S 7$, Theorem 
\ref{N-LG}) \vp 

{\bf Theorem 1} {\bf (Local-Global Principle)} Let $k$ be a commutative ring with identity
and $R$  an associative $k$-algebra such that $R$ is finite 
as a left $k$-module.
If $\alpha(X)\in {\G}(2n,R[X], \LMD[X])$, $\alpha(0)={\I}_n$ and 
$$\alpha_{\m}(X)\in {\E}(2n,R_{\m}[X], \LMD_{\mf m}[X])$$ for every maximal 
ideal $\m \in {\M}(k)$, then $$\alpha(X)\in {\E}(2n,R[X], \LMD[X]).$$ 
$($Note that $R_{\m}$ denotes $S^{-1}R$, where $S = k \setminus \m$.$)$ \vp

{\bf Theorem 2}
Let $k$ be a commutative ring with identity
and $R$  an associative $k$-algebra such that $R$ is finite 
as a left $k$-module. The for size 
at least 6 in the 
quadratic case and at least $2(r+3)$ in the Hermitian case:
\begin{center}
{\bf (Normality of the elementary subgroup)}\\ $\equiv$\\ {\bf (Local-Global Principle)}
\end{center} \vp 

\iffalse
\begin{enumerate}
 \item {\bf (Normality)} The elementary subgroup ${\E}(2n, R, \LMD)$ is a normal subgroup of 
the general quadratic (general Hermitian) group ${\G}(2n, R, \LMD)$.
\item {\bf (L-G Principle)}
If $\alpha(X)\in {\G}(2n,R[X], \LMD[X])$, $\alpha(0)={\I}_n$ and 
$$\alpha_{\m}(X)\in {\E}(2n,R_{\m}[X], \LMD_{\mf m}[X])$$ for every maximal 
ideal $\m \in {\M}(k)$, then $$\alpha(X)\in {\E}(2n,R[X], \LMD[X]).$$ 
$($Note that $R_{\m}$ denotes $S^{-1}R$, where $S = k \setminus \m$.$)$
\end{enumerate}
\fi 

To give a complete picture about the  ${\k}$-functors 
we shall shorty discuss the progress in the stabilization problem for ${\k}$-functors.  The study of this 
problem first time appeared in the work of 
Bass--Milnor--Serre, and then we see in the work by  A. Bak, M. Stein, L.N. Vaserstein, and others 
for the symplectic, orthogonal and general quadratic groups. For details {\it cf.}~\cite{Bak1}, \cite{ST}, 
\cite{V}, \cite{V2}, and \cite{V3}. In 1998, R.A. Rao and W. van der Kallen 
studied this problem for the linear 
groups over an affine algebra in \cite{RV}. 
The result settled for the general quadratic and the general Hermitian groups by A. Bak, G. Tang and V. Petrov in 
\cite{Bak2} and \cite{Bak3}. The result by Bak---Petrov---Tang has been
improved by Sergei Sinchuk, ({\it cf.}~\cite{SS}). 
It has been observed that 
over a regular affine algebra Vaserstein's bounds for the stabilization can be improved for the 
transvection subgroup of full automorphism group of projective, and 
symplectic modules. But cannot be improved for the orthogonal case in general. For details {\it cf.}
\cite{BR}, \cite{BRS}. We refer recent breakthrough result by 
J. Fasel, R.A. Rao and R.G. Swan (\cite{FRS}, Corollary 7.7).
%showed that if  $R$ is a smooth algebra of dimension $d\ge 3$ over an algebraically closed
%field $k$, then, with the assumption  $d!k = k$, ${\SL}_d (R)\cap {\E}_{d+1} (R) = {\E}_d (R)$ (\cite{FRS}, Corollary 7.7). 
A very recent result of Weibo Yu gives a similar bound for the odd unitary groups, ({\it cf.}~\cite{We}). 
In this paper we don't prove any new result in this direction. 

%In the future paper, 
%the author expects to establish an analogous local-global principle for the transvection 
%subgroup of the full automorphism group in the case of the general quadratic group. As an 
%application it will follow that in the case of the general quadratic group 
%one can deduce stabilization bound for the module case from the known bound of the 
%free case. The author believes that the same technique will be applicable for any other 
%classical-like groups. 

Though the study of stability for ${\k}$-functors started in mid-sixties, the structure of 
${\k}$-group below the level of stable range was not much studied. 
In 1991, A. Bak  showed that 
the group ${\GL}(n, R)/{\E}(n, R)$ is nilpotent-by-abelian for $n\ge 3$; ({\it cf.}~\cite{Bak}). 
In \cite{HR}, R. Hazrat proved the similar result for the general quadratic groups over module finite rings. 
The paper of Hazrat and Vavilov \cite{HV} redoes this for classical Chevalley groups (that is types A, C, and D)  
and then extends it further to the exceptional Chevalley groups (that is types E, F, and G). 
They have shown the following: Let $\Phi$ be a reduced
irreducible root system of rank $\geq 2$ and $R$ be a commutative ring
such that its Bass--Serre dimension $\delta(R)$ is finite. Then for any   
Chevalley group ${\G}(\Phi, R)$ of type $\Phi$ over $R$ the quotient 
${\G}(\Phi, R)/{\E}(\Phi, R)$ is nilpotent-by-abelian. In particular, 
${\k}(\Phi, R)$ is nilpotent of class at most $\delta(R) + 1$. 
They use the localization-completion method of A. Bak in \cite{Bak}. 
In \cite{BBR}, the author with Bak and Rao gave a uniform proof for 
the transvection subgroup of 
full automorphism group of projective, 
symplectic and orthogonal modules of global rank at least 1 and 
local rank at least 3. Our method of proof shows that for classical groups the 
localization part suffices. Recently, in ({\it cf.}~\cite{BHV}) Bak, Vavilov and 
Hazrat proved the relative case for the unitary and Chevalley groups. 
But, to my best knowledge, so far there is no definite result for the general Hermitian groups. 
I observe that using the above local-global principle,
arguing as in \cite{BBR}, it follows that the unstable ${\k}$ of general Hermitian 
group is nilpotent-by-abelian. We follow the line of Theorem 4.1 in \cite{BBR}. 
More precisely, we prove ($\$ 8$, Theorem \ref{nil}) \vp 

{\bf Theorem 3} For the general Hermitian group of large size  over a commutative ring $R$ with identity, 
the quotient 
group $\frac{{\SH}(2n, R, a_1,\ldots, a_r)}{{\EH}(2n, R, a_1,\ldots, a_r)}$ is 
nilpotent for $n\ge r+3$. \vp 

We conclude with a brief description of the organization 
of the rest of the paper.
Section 1 of the paper serves as an introduction.
Section 2 recalls form rings,  section 3 general quadratic groups 
over form rings and their elementary subgroups, section 4 general Hermitian groups and their 
elementary subgroups, section 5 provides preliminary results regarding the groups above, 
section 6 the local-global principle for the elementary subgroup of the general quadratic and 
general Hermitian group, section 7 the equivalence of normality of the elementary subgroup 
and the local-global principle for the elementary subgroup, and section 8 the nilpotent by 
abelian structure of non-stable ${\k}$ of the general Hermitian group.

\section{Form Rings} 
{\bf Definition:} 
Let us first recall the concept of $\LMD$-quadratic forms introduced 
by A. Bak in his Ph.D. thesis ({\it cf.}~\cite{Bak1}) in order to overcome the difficulties that arise for the 
characteristic 2 cases. 

Let $R$ be an (not necessarily commutative) associative ring with identity, and 
with involution $-:R\ra R$, $a\mapsto \ol{a}$. 
Let $\lambda\in C(R)$ = center of $R$ be an element with the property $\lmd \ol{\lmd}=1$. 
We define additive subgroups of $R$
$$\LMD_{\tn{max}}=\{a\in R\,|\, a=-\lmd\ol{a}\} \,\,\,\, \& \,\,\,\,
\LMD_{\tn{min}}=\{a-\lmd\ol{a}\,|\, a\in R\}.$$ 

One checks that $\LMD_{\tn{max}}$ and $\LMD_{\tn{min}}$
are closed under the conjugation operation $a\mapsto \ol{x}ax$ for any
$x\in R$. A $\lmd$-form parameter on $R$ is an additive subgroup $\LMD$ of $R$
such that $ \LMD_{\tn{min}}\subseteq \LMD\subseteq\LMD_{\tn{max}}$, and 
$\ol{x}\LMD x\subseteq \LMD$ for all $x\in R$. A pair $(R,\LMD)$ is called a {\it form ring}. \vp 

{\bf Examples:} 
\begin{enumerate} 
\item $ \LMD_{\tn{min}} = 0 \Leftrightarrow \lmd=1$, and involution is trivial.  
In particular, $\LMD=0 \Leftrightarrow \lmd=1$, involution is trivial, and $R$ is commutative.
\item  If $R$ is a commutative integral domain, and involution is trivial. {\it i.e.}, then  
$\lmd^2=1 \Leftrightarrow \lmd=\pm 1$. If 
$\lmd=1$ and char$R\ne 2$, then $\LMD_{\tn{max}}=0$, and so $0$ is the only form parameter. If 
$\lmd=-1$ and char$R\ne 2$, then $\LMD$ contains $2R$, and closed under multiplication by squares. 
If  $R$ is a field, then we get $\LMD = R$. 
If $R$ is a $\mbb{Z}$, then we get $\LMD = 2\mbb{Z}$ and $\mbb{Z}$. 
If char$R = 2$, then $R^2$ is a subring of $R$, and $\LMD$ = $R^2$-submodules of $R$.
\item The ring of $n\times n$ matrices $({\EM}(n, R), \LMD_n)$ is a form ring. 
\end{enumerate} \vp

%{\bf Change of Bases:} 
%\begin{enumerate}
%\item Polynomial Extension: Let $(R[X], \LMD[X])$ be the form ring induced from $(R, \LMD)$. 
%\item Localization: Let $(S^{-1}R[X], S^{-1}\LMD[X])$ be the form ring induced from $(R, \LMD)$ after localization. 
%\end{enumerate} \vp 

{\bf Remark:} Earlier version of  $\lmd$-form parameter  
is due to K. McCrimmon which plays an important role in his 
classification theory of Jordan Algebras. He defined for the wider class of 
alternative rings (not just associative rings), but for associative rings it 
is a special case of Bak's concept. (For details, {\it cf.} N. Jacobson; Lectures on 
Quadratic Jordan Algebras, TIFR, Bombay 1969).  The excellent work of 
Hazrat-Vavilov in \cite{HV1} is a  very good source to understand the historical motivation 
behind the concept of form rings. And, an excellent source to understand the theory of
form rings is the book \cite{HO} by A.J. Hahn and O.T. O'Meara.

\section{General Quadratic Group}
Let $V$ be a right $R$-module and 
${\GL}(V)$ the group of all $R$-linear automorphisms of $V$. 
A map $f: V\times V \ra R$ is called 
{\it sesqulinear form} if 
$f(ua,vb)=\ol{a}f(u,v)b$ for all $u,v\in V$, \,\,$a,b\in R$. 
We define $\LMD$-{\it quadratic form} $q$ on $V$, and associated $\lmd$-{\it Hermitian form} and  as follows:  
$$q : V \ra R/\LMD, \,\,{\rm ~given~ by~} \,\, q(v)=f(v,v) + \LMD,\,\,{\rm ~and~}$$  
$$h: V\times V\ra R; \,\,{\rm ~given~ by~}\,\, h(u,v)=f(u,v) + \lmd \ol{f(v,u)}.$$
A Quadratic Module over $(R, \LMD)$ is a triple $(V,h,q)$.  \vp 

{\bf Definition:} 
``Bak's Unitary Groups'' or ``The Unitary Group of a Quadratic Module'' or 
``General Quadratic Group'' ${\GQ}(V,q,h)$ is defined as follows: 
$${\GQ}(V,q,h) \,\, =\,\,
\{\alpha \in {\GL}(V) \,\,|\,\, h(\alpha u, \alpha v)=h(u,v), \,\, q(\alpha u)=q(v)\}.$$

{\bf Examples:} \underline{Traditional Classical Groups}
\begin{enumerate} 
\item By taking $\LMD = \LMD_{\tn{max}} = R$, $\lmd=-1$, and trivial involution we get
symplectic group ${\GQ}(2n, R, \LMD) = {\Sp}(2n, R)$.
\item By taking $\LMD = \LMD_{\tn{min}} = 0$, $\lmd=1$, trivial involution we get
quadratic or orthogonal group ${\GQ}(2n, R, \LMD) = {\O}(2n, R)$.
\item For general linear group, let $R^{o}$ be the ring opposite to $R$, 
and $R^e = R\oplus R^{o}$. Define involution as follows: $(x, y^o) \mapsto (y, x^o)$. 
Let $\lmd = (1,1^o)$ and $\LMD = \{(x,-x^o) \,\,|\,\,x\in R\}$. Then identify 
${\GQ}(2n, R^e, \LMD) = \{(g, g^{-1})\,\,|\,\, g\in {\GL}(n,R)\}$ with ${\GL}(n,R)$. 
\end{enumerate}

%{\bf Elementary Quadratic Transvections:} \vp 

{\bf Free Case:} Let $V$ be a free right $R$-module of rank $2n$ 
with ordered basis $e_1, e_2, \ldots, e_n, e_{-n},\ldots, e_{-2}, e_{-1}$. 
Consider the sesqulinear form 
$f : V\times V \lra R$, defined by $f(u,v)=\ol{u}_1v_{-1}+\cdots+\ol{u}_nv_{-n}$. 
Let 
$h$ be Hermitian form, and 
$q$ be the $\LMD$-quadratic form defined by $f$. So, we have 
$$h(u,v) = \ol{u}_1v_{-1}+\cdots+\ol{u}_nv_{-n} +\lmd \ol{u}_{-n}v_n+\cdots+\lmd\ol{u}_{-1}v_1,$$ 
$$q(u) = \LMD + \ol{u}_1u_{-1}+\cdots+\ol{u}_n u_{-n}.$$

Using this basis we can identify ${\GQ}(V, h, q)$ with a subgroup of ${\GL}(2n, R)$ of rank $2n$. 
We denote this subgroup by ${\GQ}(2n, R, \LMD)$. 

By fixing a basis $e_1,e_2,\ldots, e_n, e_{-1}, e_{-2}, \ldots, e_{-n}$, 
we define the form
$$\psi_n= \begin{pmatrix} 0 & \lmd {\rm I}_n \\ {\rm I}_n &0\end{pmatrix}$$ 
Hence, 
${\GQ}(2n, R,\LMD)$ = $\{\sigma\in {\GL}(2n, R, \LMD)\,|\, \ol{\sigma}\psi_n
\sigma=\psi_n\}.$ \vp 

For $\sigma=\begin{pmatrix} \alpha & \beta \\ \gamma &
\delta \end{pmatrix}\in {\GL}(2n, R,\LMD)$, one can show that 
$\sigma\in {\GQ}_{2n}(R,\LMD)$ ($\alpha, \beta, \gamma, \delta$ are $n\times n$ block 
matrices) if and only if $\ol{\gamma} \alpha,\ol{\delta}\beta\in \LMD$. 
For more details see (\cite{Bak1}, 3.1 and 3.4). \vp 

A typical element in ${\GQ}(2n, R,\LMD)$  
is denoted by a $2n\times 2n$ matrix $\begin{pmatrix} \alpha & \beta \\ 
\gamma &\delta \end{pmatrix}$, where $\alpha, \beta, \gamma, \delta$ are
$n\times n$ block matrices. 

There is a standard embedding, 
${\GQ}(2n, R,\LMD)\ra {\GQ}(2n+2, R,\LMD)$, given by 
$$\begin{pmatrix} \alpha & \beta \\ 
\gamma &\delta \end{pmatrix} \mapsto \begin{pmatrix}\alpha & 0 & \beta & 0\\
0 & 1 & 0 & 0\\ \gamma & 0 & \delta & 0\\0 & 0 & 0 & 1\end{pmatrix}$$
called the {\it stabilization} map. This allows us to identify ${\GQ}(2n, R,\LMD)$ 
with a subgroup in ${\GQ}(2n+2, R,\LMD)$. \vp 

{\bf  Elementary Quadratic Matrices:}
Let $\rho$ be the permutation, defined by $\rho(i)=n+i$ for $i=1,\ldots,n$. 
%Let $e_i$ denote the column vector with $1$ in the $i$-th position and 0's elsewhere. 
Let $e_{ij}$ be the matrix with $1$ in the $ij$-th position and 0's
elsewhere. For $a\in R$, and  $1\le i, j\le n$, we define 
$$q\eps_{ij}(a)={\rm I}_{2n}+ae_{ij}-\ol{a}e_{\rho(j)\rho(i)} \,\,\tn{ ~~~~for } i\ne j,$$
$$qr_{ij}(a)=\begin{cases} {\rm I}_{2n}+ae_{i\rho(j)}-\lmd\ol{a}e_{j\rho(i)} \,\,\tn{ for }  i\ne j\\
{\rm I}_{2n}+ae_{\rho(i)j}\,\,\,\,\,\,\,\,\,\,\,\ \,\,\,\,\,\,\,\,\,\,\,\,\, \,\,\,\tn{~~ for } 
i=j,\end{cases}$$ 
$$ql_{ij}(a)=\begin{cases}
{\rm I}_{2n}+ae_{\rho(i)j}-\ol{\lmd}\ol{a}e_{\rho(j)i} \,\,\tn{ for } i\ne j\\
{\rm I}_{2n}+ ae_{\rho(i)j} \,\,\,\,\,\,\,\,\,\,\,\ \,\,\,\,\,\,\,\,\,\,\,\,\, \,\,\,\tn{ ~~for } 
i=j.\end{cases}$$
(Note that for the second and third type of elementary matrices, if $i=j$, 
then we get $a=-\lmd\ol{a}$, and hence it forces that $a\in \LMD_{\rm max} (R)$. 
One checks that these
above matrices belong to ${\GQ}(2n, R,\LMD)$; {\it cf.}~\cite{Bak1}.) \vp 

{\bf n-th Elementary Quadratic Group
${\EQ}(2n, R,\LMD)$:} \tn{The subgroup generated by $q\eps_{ij}(a)$, 
$qr_{ij}(a)$ and $ql_{ij}(a)$, for $a\in R$ and $1\le i,j\le n$.} \vp 

It is clear that the stabilization map takes 
generators of ${\EQ}(2n, R,\LMD)$ to the generators of 
${\EQ}(2(n+1), R, \LMD)$. 

{\bf Commutator Relations:} 
There are standard formulas for the commutators between quadratic 
elementary matrices. For details we refer \cite{Bak1} (Lemma 3.16),
and \cite{HR} ($\S$ 2).   
In later sections we shall repeatedly use those relations. \vp

\section{Hermitian Group}

We assume that $\LMD$ is a $\lmd$ form parameter on $R$. 
For a matrix $M=(m_{ij})$ over $R$ we  
define $\ol{M}=(\ol{m}_{ij})^t$. 
For $a_1,\dots,a_n\in \LMD$ and $n>r$ let 
$$A_1=\begin{pmatrix} a_1 & 0 & 0 & \cdots & 0\\
0 & a_2 & 0 & \cdots & 0\\
\cdots & \cdots & \cdots & \cdots & \cdots\\
0 & \cdots & 0 & a_{r-1} & 0\\
0 &  \cdots & 0 &  0 & a_r
\end{pmatrix}=[a_1,\ldots,a_r]$$ 
denote the diagonal matrix whose $ii$-th diagonal coefficient is $a_i$. 
Let $A=A_1\perp {\rm I}_{n-r}$. We define the form
$$\psi^h_n= \begin{pmatrix} A_1 & \lmd {\rm I}_n \\ {\rm I}_n &0\end{pmatrix}.$$ 

{\bf Definition:}  {\bf General Hermitian Group} 
\tn{of the elements $a_1,\ldots,a_r$ is defined as follows: 
${\GH}(2n, R, a_1,\ldots ,a_r, \LMD)$: The group generated by the all 
non-singular $2n\times 2n$ matrices} $$\{\sigma\in 
{\GL}(2n, R)\,|\, \ol{\sigma}\psi^h_n \sigma=\psi^h_n\}.$$

As before, there is an obvious embedding 
$${\GH}(2n, R, a_1,\ldots,a_r, \LMD) \inj {\GH}(2n+2, R, a_1,\ldots,a_r, \LMD).$$

To define {\bf elementary Hermitian matrices}, we need to consider the set 
$C=\{(x_1,\ldots,x_r)^t\in (R^r)^t\,|\,\underset{i=1}{\overset{r}\sum}
\ol{x}_ia_ix_i\in \LMD_{\tn{min}} (R)\}$
for $a_1,\ldots,a_r$ as above. In order to overcome the technical difficulties 
caused by the elements $a_1,\ldots,a_r$, we shall finely partition a typical
matrix $\begin{pmatrix} \alpha & \beta \\ 
\gamma &\delta \end{pmatrix}$ of ${\GH}(2n, R, a_1,\ldots,a_r, \LMD)$ into the form 
$$\begin{pmatrix} \alpha_{11} & \alpha_{12} & \beta_{11} & \beta_{12} \\
\alpha_{21} & \alpha_{22} & \beta_{21} & \beta_{22}\\
\gamma_{11} & \gamma_{12} & \delta_{11} & \delta_{12} \\
\gamma_{21} & \gamma_{22} & \delta_{21} & \delta_{22} \end{pmatrix}$$  
where $\alpha_{11},\beta_{11}, \gamma_{11}, \delta_{11}$ are $r\times r$
matrices, $\alpha_{12},\beta_{12},\gamma_{12},\delta_{12}$ are $r\times
(n-r)$ matrices, $\alpha_{21},\beta_{21},\gamma_{21},\delta_{21}$
are $(n-r)\times r$ matrices, and $\alpha_{22},\beta_{22},\gamma_{22},
\delta_{22}$ are $(n-r)\times (n-r)$ matrices. By (\cite{T}, Lemma 3.4), 
\begin{eqnarray} \label{tn1} 
\tn{  the columns of  }
\alpha_{11}-{\rm I}_r,\alpha_{12},\beta_{11},\beta_{12},\ol{\beta}_{11},
\ol{\beta}_{21}, \ol{\delta}_{11}-{\rm I}_r,\ol{\delta}_{21} \in C. 
\end{eqnarray} 
It is a straightforward check that the subgroup of \lb ${\GH}(2n, R, a_1,\ldots,a_r, \LMD)$ 
consisting of  
\begin{eqnarray*}
  \left\{ \begin{pmatrix} {\rm I}_r & 0 & 0 & 0\\
0 & \alpha_{22} & 0 & \beta_{22}\\
0 & 0 & {\rm I}_r & 0 \\
0 & \gamma_{22} & 0 & \delta_{22} \end{pmatrix}
\in {\GH}(2n, R, a_1,\ldots,a_r) \right\}
\end{eqnarray*} 
\begin{eqnarray*}
\cong {\GH}(2(n-r), R,a_1,\ldots,a_r,\LMD).
\end{eqnarray*} 

{\bf  Elementary Hermitian Matrices:}
The first three kinds of generators are taken for the most part from 
${\GQ}(2(n-r), R,\LMD)$, which is embedded, as above, as a subgroup of
${\GH}(2n, R)$ and the last two kinds are motivated by the result \eqref{tn1}
concerning the column of a matrix in ${\GH}(2n, R)$. For $a\in R$, we define 

\begin{align*}
h\eps_{ij}(a)= & {\rm I}_{2n}+ae_{ij}-\ol{a}e_{\rho(j)\rho(i)} \,\,\,\,~~\tn{~~~~~for}\,\,
 r+1\le i\le n, 1\le j\le n, i\ne j, \\
hr_{ij}(a)= & \begin{cases} 
{\rm I}_{2n}+ae_{i\rho(j)}-\lmd\ol{a}e_{j\rho(i)}  \,\,\,\,\tn{for}\,\,  r+1\le i, j\le n, i\ne j\\
{\rm I}_{2n}+ae_{i\rho(j)} \,\,\,\,\,\,\,\,\,\,\,\,\,\,\,\,\,\,\,\,\,\,\,\,\,\,\,\,\,\,\tn{~for} 
\,\, r+1\le i, j\le n, i=j,
\end{cases} \\
hl_{ij}(a)= & \begin{cases} {\rm I}_{2n}+ae_{\rho(i)j}-\ol{\lmd}\ol{a}e_{\rho(j)i} 
\,\, \,\,\tn{for} \,\,  1\le i, j\le n, i\ne j \\
 {\rm I}_{2n}+ae_{\rho(i)j}  \,\,\,\,\,\,\,\,\,\,\,\,\,\,\,\,\,\,\,\,\,\,\,\,\,\,\,\,\,\,
\tn{~for} \,\,  1\le i, j\le n, i= j. 
\end{cases}
\end{align*}
(Note that for the second and third type of elementary matrices, if $i=j$, 
then we get $a=-\lmd\ol{a}$, and hence it forces that $a\in \LMD_{\rm max} (R)$).  One checks that the 
above matrices belong to ${\GH}(2n, R,a_1,\ldots,a_r, \LMD)$; {\it cf.}~\cite{T}. 

For $\zeta=(x_1,\ldots,x_r)^t\in C$, let $\zeta_f\in R$ be such that 
$\zeta_f+\lmd\ol{\zeta}_f=\underset{i=1}{\overset{r}\sum}
\ol{x}_ia_ix_i$. (The element $\zeta_f$ is not unique in general). We define
$$hm_i(\zeta)=\begin{pmatrix} {\rm I}_r & \alpha_{12} & 0 & 0\\
0 & {\rm I}_{n-r} & 0 & 0\\
0 & -\ol{A}_1\alpha_{12} & {\rm I}_r & 0 \\
0 & \gamma_{22} & -\ol{\alpha}_{12} & {\rm I}_{n-r} \end{pmatrix}$$ for $\zeta\in C \,\,\& \,\, r+1\le i\le n$
to be the  $2n\times 2n$ matrix, where $\alpha_{12}$ is the $r\times (n-r)$ matrix with
$\zeta$ as its $(i-r)$-th column and all other column's zero, and $\gamma_{22}$ is
the $(n-r)\times (n-r)$ matrix with $\ol{\zeta}_f$ in $(i-r,i-r)$-th position and $0$'s elsewhere. 
Let $e_k$ denote the column vector of length $(n-r)$ with $1$ in the $k$-th position and $0$'s elsewhere, 
and $e_{t\,s}$ denote a $(n-r)\times (n-r)$ matrix with $1$ in the $ts$-th position and $0$'s elsewhere. 
%Thus $\alpha_{12} e_{i-r} = \zeta$. 
%$\gamma_{22}e_{(i-r)\,\,(i-r)}=\ol{\zeta}_f$. 

As above, we define 
$$hr_i(\zeta)=\begin{pmatrix} {\rm I}_r & 0 & 0 & \beta_{12}\\
0 & {\rm I}_{n-r} & -\lmd\ol{\beta}_{12} & \beta_{22}\\
0 & 0 & {\rm I}_r & -\ol{A}_1\beta_{12} \\
0 & 0 & 0 & {\rm I}_{n-r} \end{pmatrix}$$ for  
$\zeta\in C \,\,\& \,\, r+1\le i\le n$
a $2n\times 2n$ matrix, where $\beta_{12}$ 
is the $r\times (n-r)$ matrix with
$\zeta$ as its $(i-r)$-th column and all other column's zero, and $\beta_{22}$ is
the $(n-r)\times (n-r)$ matrix with 
$\lmd\ol{\zeta}_f$ in $(i-r,i-r)$-th position and $0$'s elsewhere. 
%Thus $\beta_{12}e_{i-r}=\zeta$.
%$\beta_{22}e_{(i-r)\,\,(i-r)}=\lmd \ol{\zeta}_f$.  

Note that if $\eta=e_{pq}(a)$ is an elementary generator in ${\GL}(s, R)$, 
then the matrix $({\rm I}_{n-s}\perp
\eta\perp {\rm I}_{n-s} \perp \eta^{-1})\in h\eps_{ij}(a)$. 
It has been shown in \cite{T} ($\S 5$) that each of the 
above matrices is in ${\GH}(2n, R, a_1,\ldots,a_r, \LMD)$.\vp 

{\bf Definition:} {\bf n-th Elementary Hermitian Group} 
of the elements $a_1,\ldots, a_r$;
${\EH}(2n, R,a_1,\ldots,a_r, \LMD)$: The group generated by $h\eps_{ij}(a)$, 
$hr_{ij}(a)$, $hl_{ij}(a)$, $hm_i(\zeta)$ and $hr_i(\zeta)$, for $a\in R$, $\zeta\in C$ 
and $1\le i,j\le n$.

The stabilization map takes 
generators of  ${\EH}(2n, R,a_1,\dots, a_r,\LMD)$ to the  generators of 
${\EH}(2(n+1), R, a_1,\dots, a_r,\LMD)$. \vp

{\bf Commutator Relations:} There are standard formulas for the commutators between quadratic 
elementary matrices. For details we refer  \cite{T}.  \vp

\section{Preliminaries and Notations}

{\bf Blanket Assumption:} We always assume that $2n\ge 6$ and  
$n>r$ while dealing with the Hermitian case. 
%We consider only isotropic vectors. 
We do not want to put any restriction on the
elements of $C$. Therefore we assume that $a_i\in \LMD_\tn{min}(R)$ for 
$i=1,\ldots,r$, as in that case $C=R^r$. We always assume 
$a_1=0$.
\begin{nt} \tn{In the sequel \tn{M}$(2n,R)$ will denote the set of all 
$2n\times 2n$ matrices. By $\G(2n,R, \LMD)$ we shall denote either the quadratic group 
${\GQ}(2n, R,\LMD)$ or the Hermitian group ${\GH}(2n, R,a_1,\ldots,a_r, \LMD)$ of size $2n\times 2n$. By
${\es}(2n, R, \LMD)$ we shall denote respective subgroups ${\SQ}(2n, R,\LMD)$ or 
${\SH}(2n, R,a_1,\ldots,a_r,\LMD)$ with matrices of determinant 1, 
in the case when $R$ will be commutative.  Then, by
${\E}(2n, R, \LMD)$ we shall denote the corresponding elementary subgroups
${\EQ}(2n, R,\LMD)$ and ${\EH}(2n, R,a_1,\ldots,a_r, \LMD)$. 
To treat uniformly we denote the elementary generators of ${\EQ}(2n, R,\LMD)$, 
and the first three types of 
elementary generators of ${\EH}(2n, R,\LMD)$ by $\vartheta_{ij}(\star)$, for some $\star\in R$.
To express the last two types of generators of ${\EH}(2n, R,\LMD)$ we shall use the notation 
$\vartheta_i(\star)$, where $\star$ is a column vector of length $r$ defined over the ring $R$. 
{\it i.e.}, we will have two types of elementary generators, namely 
$\vartheta_{ij}({\rm ring \,\, element})$ and $\vartheta_i({\rm column \,\, vector})$. 
Let $\LMD[X]$ denote the $\lambda$-form parameter on $R[X]$ induced from $(R,\LMD)$, {\it i.e.}, $\lmd$-form parameter 
on $R[X]$ generated by $\LMD$, {\it i.e.}, the smallest form parameter on $R[X]$ containing $\LMD$. 
Let $\LMD_s$ denote the $\lambda$-form parameter on $R_s$ induced from $(R,\LMD)$. }
\end{nt} 

For any column vector $v\in (R^{2n})^t$ we define the row vectors
$\widetilde{v}_q=\ol{v}^t\psi^q_n$ and $\widetilde{v}_h=\ol{v}^t\psi^h_n$.
\begin{de} \tn{We define the map $\tn{M}:(R^{2n})^t \times (R^{2n})^t \ra \tn{M}(2n,R)$ and 
the inner product $\langle \,,\rangle$ as follows: }
\begin{align*}
\tn{M}(v,w) & = v.\widetilde{w}_q-\ol{\lambda}\ol{w}.\widetilde{v}_q, \,\,\,\, \tn{ when } 
      {\G}(2n,R)={\GQ}(2n, R,\LMD)\\
  & = v.\widetilde{w}_h -\ol{\lambda}\ol{w}.\widetilde{v}_h, \,\,\,\,  \tn{ when } 
      {\G}(2n,R)={\GH}(2n, R,a_1,\ldots,a_r, \LMD),\\ 
\langle v,w\rangle & = \widetilde{v}_q.w, \,\,\,\, \tn{ when } 
{\G}(2n,R)={\GQ}(2n, R,\LMD) \\
  & = \widetilde{v}_h.w, \,\,\,\, \tn{ when } 
             {\G}(2n,R)={\GH}(2n, R,a_1, \ldots , a_r, \LMD).
\end{align*} 
\end{de} 

Note that the elementary generators of the groups ${\EQ}(2n, R)$ and ${\EH}(2n, R)$ are of the form 
${\rm I}_{2n}+{\rm M}(\star_1,\star_2)$ for suitable chosen standard basis vectors. 
%In fact, it is precisely the set $$\{{\rm I}_{2n}+{\rm M}(\star_1,\star_2)| \tn{ where } 
%\star_1, \star_2 \tn{ are standard basis vectors}, \langle \star_1, \star_2 \rangle = 0\}.$$  

We recall the following well known facts: 

\begin{lm} \label{nor} $($cf. \cite{Bak1}, \cite{T}$)$
The group  ${\E}(2n,R, \LMD)$ is perfect for $n\ge 3$ in the quadratic case, and 
for $n\ge r+3$ in the Hermitian case, {\it i.e.}, 
$$[{\E}(2n,R, \LMD), {\E}(2n,R, \LMD)] = {\E}(2n,R, \LMD).$$ 
\end{lm}

\begin{lm} \label{key2} {\bf (Splitting property):} 
For all elementary generators of the general quadratic group ${\GQ}(2n, R,\LMD)$ and for the first three 
types elementary generators of the Hermitian group ${\GH}(2n, R,a_1,\ldots,a_r, \LMD)$ we have: 
$$\vartheta_{ij}(x+y)=\vartheta_{ij}(x)\vartheta_{ij}(y)$$
for all $x, y \in R$.  

For the last two types of  elementary generators of Hermitian group we have the following relation: 

$$hm_i(\zeta)hm_i(\xi)=hm_i(\zeta + \xi)hl_{ii}(\bar{\zeta}_f+\bar{\xi}_f+\bar{\zeta}\bar{A}_1\xi-\ol{(\zeta+\xi)}_f),$$
$$hr_i(\zeta) hr_i(\xi)=hr_i(\zeta + \xi)hr_{ii}((\zeta+\xi)_f-\xi_f-\zeta_f-\bar{\xi}A_1\zeta).$$
\end{lm}

{\bf Proof.}  See pg. 43-44, Lemma 3.16, \cite{Bak1} for the ${\GQ}(2n, R,\LMD)$ 
and Lemma 8.2, \cite{T} for the group ${\GH}(2n, R,a_1,\ldots,a_r, \LMD)$.

\begin{lm} \label{key6} 
Let $G$ be a group, and $a_i$, $b_i \in G$, for 
$i = 1, \ldots, n$. Then  for $r_i={\underset{j=1}
{\overset{i}\Pi}} a_j$, we have
${\underset{i=1}{\overset{n}\Pi}}a_i b_i=
{\underset{i=1}{\overset{n}\Pi}}r_ib_ir_i^{-1} 
{\underset{i=1}{\overset{n}\Pi}}a_i.$
\end{lm}
\begin{nt} \tn{ By ${\G}(2n,R[X], \LMD [X], (X))$ we shall mean the group of all quadratic and Hermitian matrices
over $R[X]$ which are ${\I}_n$ modulo $(X)$.} 
\end{nt} 
\begin{lm} \label{key1} The group 
${\G}(2n,R[X],\LMD [X], (X)) \cap {\E}(2n,R[X], \LMD [X])$
is generated by the elements of the types 
$ \eps \vartheta_{ij}(\star_1)\eps^{-1}$ and $ \eps \vartheta_{i}(\star_2)\eps^{-1}$, where  $\eps \in {\E}(2n,R, \LMD)$,
$\star_1\in R[X]$, $\star_2\in ((R[X])^{2n})^t$ with both $\vartheta_{ij}(\star_1)$ and 
$\vartheta_{i}(\star_2)$ congruent to 
${\rm I}_{2n}$ modulo $(X)$.
\end{lm}

We give a proof of this Lemma for the Hermitian group. The proof for the quadratic case is similar, but easier. 

{\bf Proof of Lemma \ref{key1}}. Let $a_1(X),\ldots, a_r(X)$ be $r$ 
elements in the polynomial ring $R[X]$ with 
respect to which we are considering the Hermitian group 
${\GH}(2n, R[X], a_1(X),\ldots,a_r(X), \LMD [X])$. 

Let $\alpha(X)\!\in\!{\EH}(2n, R[X], a_1(X),\ldots,a_r(X), \LMD [X])$ be such that 
$\alpha(X)$ is congruent to ${\rm I}_{2n}$ modulo $(X)$. 
Then we can write $\alpha(X)$ as a product of elements of the 
form $\vartheta_{ij}(\star_1)$, where $\star_1$ is 
a polynomial in $R[X]$, and of the form $\vartheta_{i}(\star_2)$, 
where $\star_2$ is a column vector of length $r$ defined 
over $R[X]$. We write each $\star_1$ as a sum of a constant term and a 
polynomial which is identity modulo $(X)$. 
Hence by using the splitting property described in Lemma \ref{key2} each 
elementary generator $\vartheta_{ij}(\star_1)$ 
of first three type can be written as a product of two such elementary 
generators with the left one defined on $R$ and the right 
one defined on $R[X]$ which is congruent to ${\rm I}_{2n}$ modulo $(X)$.

For the last two types of elementary generators we write each vector $\star_2$ as a 
sum of a column vector defined over the ring $R$ and a column vector of defined over 
$R[X]$ which is congruent to the zero vector of length $r$ modulo $(X)$. 
In this case, as shown in Lemma \ref{key2}, we get one extra term involving 
elementary generator of the form $hl_{ii}$ or $hr_{ii}$. But that extra term is one of the 
generator of first three types. And then we can split that term again as above. 
Therefore, $\alpha(X)$ can be expressed as a product of following types of elementary generators: 
$$\vartheta_{i j}(\star_1(0))\vartheta_{i j}(X\star_1)\,\, 
{\rm with } \,\,\star_1(0)\in R \,\,{\rm and }  \,\, \vartheta_{ij}(X\star_1)={\rm I}_{2n} 
\,\,{\rm modulo }\,\, (X),$$ 
$$\vartheta_{i}(\star_2(0))\vartheta_{i}(X\star_2) \,\,{\rm with }\,\,  
\star_2(0)\in R \,\,{\rm and }\,\, \vartheta_{i}(X\star_2)={\rm I}_{2n} \,\,\rm{  modulo  } \,\, (X).$$ 

Now result follows by using the identity described in Lemma \ref{key6}.  \hb

\section{Suslin's Local-Global Principle} 

In his remarkable thesis ({\it cf.}~\cite{Bak1}) A. Bak showed that for 
a form ring $(R, \LMD)$ the elementary subgroup 
${\EQ}(2n,R, \LMD)$ is perfect for $n\ge 3$ and hence is a normal subgroup of ${\GQ}(2n,R, \LMD)$.
As we have noted earlier, this question is related to 
Suslin's local-global principle for 
the elementary subgroup. In \cite{T}, G. Tang has shown that for 
$n\ge r+3$ the elementary Hermitian group ${\EH}(2n, R,a_1, \ldots , a_r, \LMD)$ is perfect and 
hence is a normal subgroup of ${\GH}(2n, R,a_1, \ldots , a_r, \LMD)$.
In this section we deduce an
analogous local-global principle for the elementary subgroup of the general quadratic and Hermitian 
groups, when $R$ is module finite, {\it i.e.}, finite over its center.
We use this result in $\S 5$ to prove the nilpotent property of the unstable Hermitian group
${\KH}_1$.
Furthermore, we show that
if $R$ is finite over its center, then  the normality of the elementary subgroup is
equivalent to the local-global principle. This generalizes our result in \cite{BRK}. 

The following is the key Lemma, and it tells us the reason why we need to assume that 
the size of the matrix is at least 6. In \cite{BRK},  proof is given for the 
general linear group. Arguing in similar manner by using identities of 
commutator laws result follows in the unitary and Hermitian cases. A list of commutator 
laws for elementary generators is stated in (\cite{Bak1}, pg. 43-44, Lemma 3.16) for the 
unitary groups and in (\cite{T}, pg. 237-239, Lemma 8.2) for the Hermitian groups. 
For a direct proof we refer Lemma 5, \cite{P}.

\begin{lm} \label{key3}
Suppose $\vartheta$ is an elementary generator of the general quadratic \lb 
$($Hermitian$)$ group ${\G}(2n, R[X], \LMD[X])$, $n\ge 3$. Let $\vartheta$ be 
congruent to identity modulo $(X^{2m})$, for $m>0$. Then, if we conjugate 
$\vartheta$ with an elementary generator of the general quadratic $($Hermitian$)$ 
group ${\G}(2n, R, \LMD)$, we get the resulting matrix 
is a product of elementary generators of general quadratic $($Hermitian$)$ group 
${\G}(2n, R[X], \LMD[X])$, each of which is 
congruent to identity modulo $(X^m)$. 
\end{lm}

\begin{co} \label{key4} In Lemma \ref{key3} we can take $\vartheta$ as a product of 
elementary generators of the general quadratic $($general Hermitian$)$ group ${\G}(2n, R[X], \LMD[X])$.
\end{co}

\begin{lm} \label{key5} 
Let $(R, \LMD)$ be a form ring and $v\in {\E}(2n,R, \Lambda)e_{2n}$. Let $w\in R^{2n}$ be a 
column vector such that $\langle v,w\rangle=0$. Then
 ${\rm I}_{2n}+\tn{M}(v,w)\in {\E}(2n,R, \LMD)$.
\end{lm}
{\bf Proof.} 
Let $v=\eps e_{2n}$, where $\eps \in {\E}(2n,R,\LMD)$. Then it follows 
that
${\rm I}_{2n}+\tn{M}(v,w)=\eps ({\rm I}_{2n}+ \tn{M}(e_{2n},w_1))\eps^{-1},$ where 
$w_1 = \eps^{-1}w$.
Since $\langle e_{2n},w_1\rangle=\langle v,w\rangle=0$, we get
$w_1^t=(w_{11},\dots,w_{1 \,n-1},0,w_{1\,n+1},\ldots,w_{1 \,2n})$.
Therefore, as $\lmd\ol{\lmd}=\ol{\lmd}\lmd=1$, 
$${\rm I}_{2n}+\tn{M}(v\!,\!w)\!\! = \!\! \begin{cases}
\underset{1\le i\le n-1}
{\underset{1\le j\le n} \Pi} \!\!\!\!\!
\,\eps \,ql_{in}(-\ol{\lmd} \ol{w}_{n+1\, i})\,q\eps_{jn}(-\ol{\lmd}\ol{w}_{1\, j})ql_{n\, n}^{-1}(*)
{\eps}^{-1} \\ 
\underset{1\le i\le n-1}
{\underset{r+1\le j\le n}{\underset{1\le k\le r} \Pi}} \!\!\!\!\!
\eps  hl_{in}(-\ol{\lmd} \ol{w}_{n+1\, i})\!h\eps_{jn}(-\ol{\lmd}\ol{w}_{1\, j})\!hm_n
(-\ol{w}_{1 \,k)})hl_{n\, n}^{-1}(*)
{\eps}\!\!\!^{-1} 
\end{cases}$$ 
(in the  quadratic and Hermitian cases respectively), \\
where $\ol{w}_{1 \,\,n+k}=
(w_{1 \,n+k},0,\ldots,0)$. 
Hence the result follows. \hb \vp 

Note that the above implication is true for any associative ring with identity. 
From now onwards we assume that $R$ is finite over its center $C(R)$. Let us recall 
\begin{lm} \label{noeth} 
Let $A$ be Noetherian ring and $s\in A$. Let $s\in A$ and $s\neq 0$. 
Then there exists a natural number 
$k$ such that the homomorphism ${\G}(A,s^kA, s^k\LMD) \ra {\G}(A_s, \LMD_s)$ 
$($induced by localization homomorphism $A \ra A_s)$ is injective.  
\end{lm}

For the proof of the above lemma we refer (\cite{HV}, Lemma 5.1). 
Also, we recall that one has any module finite ring $R$ is direct limit of 
its finitely generated subrings.  
%${\G}(R, \LMD) = \underset{\lra}\lim \, {\G}(R_i, \LMD_i)$, 
%where the limit is taken over all finitely generated subring of $R$. 
Thus, one may assume that $C(R)$ is Noetherian.

 Let $(R, \LMD)$ be a (module finite) form ring with identity. 

\begin{lm} \label{Di}
{\bf (Dilation Lemma)}  
Let $\alpha(X)\in {\G}(2n,R[X], \LMD [X])$, with $\alpha(0)={\rm I}_{2n}$. If 
$\alpha_s(X)\in {\E}(2n,R_s[X], \LMD_s [X])$, for some non-nilpotent $s\in R$, 
then $\alpha(bX) \in  {\E}(2n,R[X], \LMD [X])$, 
for $b\in s^l {\C}(R)$, and $l\gg 0$. 
\end{lm} 
\begin{re} 
$($In the above Lemma we actually mean there exists some $\beta(X)\in {\E}(2n,R[X], \LMD [X])$ such that
$\beta(0)={\rm I}_{2n}$ and $\beta_s(X)=\alpha(bX)$.$)$ 
%But, since there is no 
%ambiguity, for simplicity we are using the notation $\alpha(bX)$ instead of 
%$\beta_s(X)).$
\end{re} 
{\bf Proof.}  Given that $\alpha_s(X)\in {\E}(2n,R_s[X], \LMD_s [X])$.  
Since $\alpha(0)={\rm I_{2n}}$, using Lemma \ref{key1}
we can write $\alpha_s(X)$ as a product of the matrices of the form
$ \eps \vartheta_{ij}(\star_1)\eps^{-1}$ and $ \eps \vartheta_{i}(\star_2)\eps^{-1}$, 
where  $\eps \in {\E}(2n,R_s, \LMD_s)$,
$\star_1\in R_s[X]$, $\star_2\in ((R_s[X])^{2n})^t$ with both $\vartheta_{ij}(\star_1)$ and 
$\vartheta_{i}(\star_2)$ congruent to 
${\rm I}_{2n}$ modulo $(X)$. Applying the homomorphism $X\mapsto XT^d$, where $d\gg 0$, 
from the polynomial ring $R[X]$ to the polynomial ring $R[X,T]$,  
we look on $\alpha(XT^d)$. 
Note that $R_s[X,T]\cong (R_s[X])[T]$. As $C(R)$ is Noetherian, it follows from Lemma \ref{noeth} 
and Corollary \ref{key4} that over the ring 
$(R_s[X])[T]$  we can write 
$\alpha_s(XT^d)$ as a product of elementary generators of general quadratic (Hermitian) group 
such that each of those elementary generator is congruent to 
identity modulo $(T)$. Let $l$ be the maximum of the powers occurring in
the denominators of those elementary generators. Again, as $C(R)$ is Noetherian, 
by applying the homomorphism $T\mapsto s^mT$, for 
$m\ge l$, it follows from Lemma \ref{noeth} that 
over the ring 
$R[X,T]$  we can write 
$\alpha_s(XT^d)$ as a product of elementary generators of general quadratic group 
such that each of those elementary generator is congruent to identity modulo $(T)$, for some $b\in (s^l)C(R)$, 
{\it i.e.}, we get there exists some $\beta(X,T)\in {\E}(2n,R[X,T], \LMD [X,T])$ such that
$\beta(0,0)={\rm I}_{2n}$ and $\beta_s(X,T)=\alpha(bXT^d)$.
Finally, the result follows by putting $T=1$. \hb

\begin{tr} \label{LG}
{\bf (Local-Global Principle)} \\
If $\alpha(X)\in {\G}(2n,R[X], \LMD[X])$, $\alpha(0)={\I}_n$ and 
$$\alpha_{\m}(X)\in {\E}(n,R_{\m}[X], \LMD_{\mf m}[X]),$$ for every maximal ideal $\m \in \M \, 
(C(R))$, then $$\alpha(X)\in {\E}(2n,R[X], \LMD[X]).$$ $($Note that $R_{\m}$ denotes $S^{-1}R$, where 
$S = C(R) \setminus \m$.$)$
\end{tr}
{\bf Proof.} Since $\alpha_{\m}(X)\in {\E}(2n,R_{\m}[X], \LMD_{\mf m}[X])$, for all
$\m \in {\M}(C(R))$, for each $\m$ there exists $s\in C(R) \setminus \m$ such that 
$\alpha_s(X)\in  {\E}(n,R_s[X], \LMD_s[X]).$
Using Noetherian property we can consider 
a finite cover of $C(R)$, say $s_1+\cdots+ s_r=1$.  
Let $\theta(X,T)=
\alpha_s(X+T)\alpha_s(T)^{-1}$. 
Then $$\theta(X,T)\in {\E}(2n,(R_s[T])[X], \LMD_s[T][X])$$ and $\theta(0,T)={\I}_n$. 
By Dilation Lemma, applied with base ring $R[T]$, there exists $\beta(X)\in  {\E}(2n,R[X,T], \LMD[X,T])$ such that 
 $$\beta_{s}(X)= \theta(bX,T).$$ 
%$$\theta(bX,T)\in {\E}(2n,R[X,T], \LMD[X,T]) \tn{ for some } b\in (s^l)C(R),\, l\gg 0.
%\hfill\tn{ \,\, \,\,\,}(A)$$ \vp 

Since for $l\gg 0$, the ideal $\langle s_1^l,\ldots, s_r^l\rangle = R$, we chose 
 $b_1,b_2,\dots,b_r\in C(R)$, with $b_i\in (s^l)C(R),\, l\gg 0$
such that (A) holds and $b_1+\cdots+b_r=1$. 
Then there exists $\beta^i(X)\in  {\E}(2n,R[X,T], \LMD[X,T])$ such that 
 $\beta^i_{s_i}(X)= \theta(b_iX,T)$. 
Therefore, $$\underset{i=1}{\overset{r}\Pi}\beta^i(X)\in {\E}(2n,R[X,T], \LMD[X,T]).$$
 But, $$\alpha_{s_1\cdots s_r}(X)=
\left(\underset{i=1}{\overset{r-1}\Pi} \theta_{s_1\cdots \hat{s_i}\cdots s_r}(b_iX,T)
{\mid}_{T=b_{i+1}X+\cdots +b_rX}\right) \theta_{s_1\cdots \cdots s_{r-1}}(b_rX,0).$$  
Since $\alpha(0)={\I}_n$, and as a consequence of the Lemma \ref{noeth} it follows that 
the map  ${\E}(R,s^kR, s^k\LMD) \ra {\E}(R_s, \LMD_s)$ in injective, we conclude 
$\alpha(X)\in {\E}(2n,R[X], \LMD[X])$.  \hb

\iffalse
Let $b_1,b_2,\dots,b_r\in C(R)$, with $b_i\in (s^l)C(R),\, l\gg 0$
be such that (A) holds and $b_1+\cdots+b_r=1$. 
Then $\theta(b_iX,T)\in {\E}(2n,R[X,T], \LMD[X,T])$ and hence 
$\underset{i=1}{\overset{r}\Pi}\theta(b_iX,T)\in {\E}(2n,R[X,T], \LMD[X,T])$. But, 
$$\alpha(X)=\left(\underset{i=1}{\overset{r-1}\Pi} \theta(b_iX,T)
{\mid}_{T=b_{i+1}X+\cdots +b_rX}\right) \theta(b_rX,0).$$  Since 
$\alpha(0)={\I}_n$, it follows that $\alpha(X)\in {\E}(2n,R[X], \LMD[X])$.  

(For simplicity we have written the above expression. Actually we are 
considering the equality in the localization $R_{b_1.\ldots.b_r}$ and 
using Lemma \ref{noeth} repeatedly.) \hb
\fi

\section{Equivalence of Normality and Local-Global Principle} 

Next we are going to show that if $k$ is a commutative ring with identity
and $R$ is an associative $k$-algebra such that $R$ is finite 
as a left $k$-module, then the normality criterion of elementary 
subgroup is equivalent to Suslin's local-global principle for 
above two classical groups. (Remark: One can also consider $R$ as 
a right $k$-algebra.)

One of the crucial ingredients in the proof of the above theorem is the 
following result which states that the group ${\E}$ acts transitively on 
unimodular vectors. The precise statement of the fact is the following: 

\begin{de} \tn{A vector $(v_1,\ldots, v_{2n})\in R^{2n}$ is said to be unimodular if 
there exists another vector  $(u_1,\ldots, u_{2n})\in R^{2n}$ such that 
$\sum_{i=1}^{2n} v_i u_i=1$. } 

\tn{The set of all unimodular vector in $R^{2n}$ is denoted by ${\Um}(2n, R)$.} 
\end{de}

\begin{tr} \label{swan}
Let $R$ be a semilocal ring $($not necessarily commutative$)$ with involution and 
$v=(v_1,\ldots,v_{2n})^t$ be a  unimodular and isotropic vector in $R^{2n}$.
Then $v\in {\E}(2n, R)e_{2n}$ for $n\ge 2$, {\it i.e.}, ${\E}(2n,R)$ acts
transitively on the set of isotropic vectors in ${\Um}(2n, R)$. 
\end{tr}

\iffalse
\begin{co} \label{swan2}
If $R$ is a semilocal ring with involution, then one can write ${\G}(2n,R)=
{\E}(2n,R){\G}(2,R)$ for $n\ge 2$. 
%In particular, in the commutative case ${\es}(n,R)=
%{\E}(n,R){\es}(2,R)$. 
\end{co} 
\fi 
Let us first recall some known facts before we give a proof of the theorem. 

\begin{de} \tn{An associative ring $R$ is said to be {\bf semilocal} if 
$R/\tn{rad}(R)$ is Artinian semisimple.}
\end{de} 

We recall the following three lemmas. 
\begin{lm} \label{HB}
\tn{(H. Bass)} 
Let $A$ be an associative $B$-algebra such that $A$ is finite as a left
$B$-module and $B$ be a commutative local ring with identity. 
Then $A$ is semilocal. 
\end{lm} 
{\bf Proof.} Since $B$ is local, $B/\tn{rad} (B)$ is a division ring by 
definition. That implies $A/\tn{rad} (A)$ is a finite module over the division 
ring $B/\tn{rad}(B)$ and hence is a finitely generated vector space. Thus 
$A/\tn{rad} (A)$ Artinian as $B/\tn{rad}(B)$ module and hence  
$A/\tn{rad}(A)$ Artinian as $A/\tn{rad}(A)$ module, so 
it is an Artinian ring.

It is known that an artin ring is semisimple if its radical
 is trivial. Thus $A/\tn{rad}(A)$ is semisimple, as 
$\tn{rad}(A/\tn{rad}(A))=0$. 
Hence $A/\tn{rad}(A)$ Artinian semisimple. Therefore, $A$ is semilocal by 
definition. \hb

\begin{lm} \tn{(H. Bass) (\cite{B}, Lemma 4.3.26)} \label{B}
Let $R$ be a semilocal ring $($may not be commutative$)$,  and let 
$I$ be a left ideal of $R$. Let $a$ in $R$
be such that $Ra+I=R$. Then the coset $a+I=\{a+x \,|\,x\in I\}$
contains a unit of $R$. 
\end{lm}
{\bf Proof.} We give a proof due to R.G. Swan. 
We can factor out the radical and assume that $R$ is semisimple 
Artinian. Let $I=(Ra\cap I)\oplus I'$. Replacing $I$ by $I'$ we can assume 
that $R=Ra\oplus I$. Let $f:R\ra Ra$ by $r\mapsto ra$, for $r\in R$. Therefore, 
we get an split exact sequence $0\lra J\lra R \stk{f}\lra Ra\lra 0$, for some 
ideal $J$ in $R$ which gives us a map $g:R\ra J$ such that $R \stk{(f,g)}\lra 
Ra\oplus J$ is an isomorphism. Since $Ra\oplus J\cong R\cong Ra\oplus I$ 
cancellation (using Jordon-H\"older or Krull-Schmidt) shows that $J\cong I$. 
If $h:R\cong J\cong I$, then $R\stk{(f,g)}\lra Ra\oplus I\cong R$ is an \
isomorphism sending $1$ to $(a,i)$ to $a+i$, where $i=h(1)$. Hence it follows 
that $a+i$ is a unit. \hb \vp 
\begin{lm} \label{swan4} Let $R$ be a semisimple Artinian ring and $I$ be a 
left ideal of $R$. Let $J=Ra+I$. Write $J=Re$, where $e$ is an idempotent 
$($possible since $J$ is projective. For detail {\it cf.} \cite{BK} {\rm Theorem 4.2.7}$)$. 
Then there is an element $i\in I$ such that $a+i=ue$, where $u$ is a unit 
in $R$.  
\end{lm} 
{\bf Proof.} Since $R=J+R(1-e)=Ra+I+R(1-e)$, using Lemma \ref{B} 
we can find a unit $u=a+i+x(1-e)$ in $R$, for some $x\in R$. Since $a+i\in Re$, it follows that 
$ue=a+i$. \hb  

\begin{co} \label{swan5} 
Let $R$ be a semisimple Artinian ring and $(a_1,\dots,a_n)^t$ be a 
column vector over $R$, where $n\ge 2$. Let $\Sigma Ra_i=Re$, where $e$ is 
an idempotent. Then there exists $\eps\in {\E}_n(R)$ such that 
$\eps (a_1,\dots,a_n)^t=(0,\dots,0,e)^t$. 
\end{co} 
{\bf Proof.} By Lemma \ref{swan4} we can write 
$ue=\Sigma_{i=1}^{n-1} b_ia_i+a_n$, 
where $u$ is a unit. Therefore, applying an elementary transformation we can 
assume that $a_n=ue$. Multiplying from the left by 
$({\rm I}_{n-2} \perp u\perp u^{-1})$ we can make $a_n=e$. Since all $a_i$ are 
left multiple of $e$, further elementary transformations reduce our vector 
to the required form. \hb \vp 

The following observation will be needed to do the case $2n=4$. 

\begin{lm} \label{swan6} Let $R$ be a semisimple Artinian ring and $e$ be 
an idempotent. Let $f=1-e$, and $b$ be an element of $R$. If 
$fRb\subseteq Re$, then we have $b\in Re$. 
\end{lm}
{\bf Proof.} Since $R$ is a product of simple rings, it will suffice to do 
the case in which $R$ is simple. If $e=1$, we are done. Otherwise $RfR$ is a 
non-zero two sided ideal, and hence $RfR=R$. Since $Rb=RfRb\subseteq Re$, 
we have $b\in Re$. \hb 

\begin{lm} \label{swan7} Let $R$ be a semisimple Artinian ring and let
$-:R\ra R$ be a $\lmd$-involution on $R$. Let 
$( x \,\, y)^t $ be a unimodular row of length $2n$,
where $2n\ge 4$. Then there exists an element $\eps\in 
{\E}(2n,R)$ such that $\eps( x \,\, y)^t=
(x' \,\, y')^t$, where $x_1'$ is a unit
in $R$.
\end{lm} 
{\bf Proof.} Let $x=(x_1,\ldots,x_n)^t$ and 
$b=(y_1,\ldots,y_n)^t$. We claim that there exists $\eps\in {\E}(2n,R)$ 
such that $\eps( x \,\, y)^t =
(x' \,\, y')$, where $x'$ is a unit in $R$.
Among all $( x' \,\, y')^t$ of this form, choose one
for which the ideal $I=\Sigma Rx_i'$ is maximal. Replacing the original 
$(x \,\, y )^t$ by $( x' \,\, y')^t$ we can assume that $I=\Sigma Rx_i$
is maximal among such ideals. Write $I=Re$, where $e$ is an idempotent in
$R$. By Corollary \ref{swan5} we can find an element $\eta\in {\E}_n(R)$ 
such that $\eta x =(0,0,\dots,e)^t$. 
So we can modify $x$ by my elementary generators of the form $q\eps_{ij}(\star)$ or $h\eps_{ij}(\star)$ and 
hence we assume that $x=(0,0,\dots,e)^t$.
We claim that $y_i\in Re$ for all $i\ge 1$.

First we consider the case $2n\ge 6$. 
Assume $y_1\notin I$, but $y_i\in I$ for all $i\ge 2$. If 
we apply $q\eps_{1n}(1)$ in the quadratic case
then this replaces $y_{n}$ to $y_{n}-y_1$ but 
not changes $e$ and $y_1$. On the other hand for the Hermitian case we do not
have the generator $q\eps_{1n}(1)$. But if we apply $hm_n(1,\ldots,1)$, then
it changes $y_2$ but does not changes $e$ and $b_1$. 
Therefore, in both the cases we can therefore assume that some 
$y_i$ with $i>1$ is not in $I$. (Here recall that we have put no restriction 
on $C$, {\it i.e.}, for us $C=R^r$). 
Apply $qr_{i i}(1)$ with $2\le i\le n$ in the quadratic case. 
This changes $x_i=0$ (for $i>1$) to $y_i$ while $x_n=e$ is
preserved. The ideal generated by the entries of $x$ now contains $Re+Ry_i$, 
which is larger than $I$, a contradiction, as $I$ is maximal. 
In the Hermitian case if we apply suitable $hr_i(1,\ldots,1)$ then also we 
see that the ideal generated by the entries of $x$ now contains $Re+Ry_i$, 
hence a contradiction. 

If $2n=4$, we can argue as follows. Let $f=1-e$. Let us assume that 
$y_1\ne I$ as above. Then by Lemma \ref{swan6} it will follow that we can find
some $s\in R$ such that $fsy_1\ne Re$. 
First consider the quadratic case. 
Applying $qr_{21}(fs)$ replaces $x_2=e$ 
by $c=e+fsy_1$. As $ec=e$, $I=Re\subset Rc$. Also, $fc=fsy_1\in Rc$ but 
$fc\notin I$. Hence $I \subsetneq Rc$, a contradiction. We can get the similar 
contradiction for $y_2$  by applying $qr_{22}(fs)$. In the Hermitian case,
apply $hr_1(1)$ to get the contradiction for $y_1$. Now
note that in this $r=1$ as we have assume $r<n$. Hence 
we can apply $qr_{22}(fs)$ to get the contradiction.

Since all $y_i$ lie in $Re$, the left ideal generated by the all entries of 
$(x\,\,y)^t$ is $Re$, but as this column 
vector is unimodular, we get $Re=R$, and therefore $e=1$. \hb \vp \\
{\bf Proof of Theorem \ref{swan}.}  Let $J$ be the Jacobson radical of
$R$. Since the left and the right Jacobson radical are same, $J$ is stable
under the involution which therefore passes to $R/J$. Let $\eps$ be as in
Lemma \ref{swan7} for the image 
$(x'\,\,y')^t$ of
$(x\,\,y)^t$. By lifting $\eps$ from 
$R/J$ to $R$ and applying it to 
$(x\,\,y)^t$ we reduce to the case where 
$x_n$ is a unit in $R$. Let $\alpha=x_n\perp x_n^{-1}$. Then applying 
$({\rm I}_{n-2}\perp  \alpha \perp {\rm I}_{n-2} \perp \alpha^{-1})$ we can assume that 
$x_n=1$.

Next applying $\Pi_{i=1}^{n-1} ql_{ni}(-y_i)$ and  
$\Pi_{i=1}^{n-1} hl_{ni}(-y_i)$ in the respective cases we
get $y_1=\cdots=y_{n-1}=0$. As 
isotropic vector remains isotropic under
elementary quadratic (Hermitian) transformation, we have 
$y_n+\lambda\ol{y}_n=0$, hence $ql_{11}(\ol{\lmd}\ol{y}_n)$ and 
$hl_{11}(\ol{\lmd}\ol{y}_n)$ are defined and 
applying it reduces $y_n$ to $0$ in both the cases. 
Now we want to make $x_i=0$ for $i=1,\ldots,n$. In the quadratic case it can
be done by applying $\Pi_{i=1}^{n-1} h\eps_{in}(-x_i)$. Note
that this transformation does not affect any $y_i$'s, as $y_i=0$. In the
Hermitian case we can make $x_{r+1}=\cdots=x_n=0$ as before applying 
$\Pi_{i=r+1}^{n-1} q\eps_{in}(-x_i)$. To make 
$x_1=\cdots=x_r=0$ we have to recall that the set $C=R^r$, {\it i.e.}, there is
no restriction on the set $C$. Hence $hr_n(-x_1,\ldots,-x_r)$ is defined and
applying it we get $x_1=\cdots=x_r=0$. Also note that other $x_i$'s and
$y_i$'s remain unchanged. Finally, applying $hl_{nn}(1)$ and then
$hr_{nn}(-1)$ we get the required vector 
$(0,\ldots,0,1)$. This completes the proof.\hb 

\begin{tr}  \label{N-LG}
Let $k$ be a commutative ring with identity
and $R$  an associative $k$-algebra such that $R$ is finite 
as a left $k$-module. Then the following are equivalent for $n\ge 3$ in the 
quadratic case and $n\ge r+3$ in the Hermitian case:
\begin{enumerate}
 \item {\bf (Normality)} ${\E}(2n, R, \LMD)$ is a normal subgroup of 
${\G}(2n, R, \LMD)$.
\item {\bf (L-G Principle)}
If $\alpha(X)\in {\G}(2n,R[X], \LMD[X])$, $\alpha(0)={\I}_n$ and 
$$\alpha_{\m}(X)\in {\E}(2n,R_{\m}[X], \LMD_{\mf m}[X])$$ for every maximal 
ideal $\m \in {\M}(k)$, then $$\alpha(X)\in {\E}(2n,R[X], \LMD[X]).$$ 
$($Note that $R_{\m}$ denotes $S^{-1}R$, where $S = k \setminus \m$.$)$
\end{enumerate}
\end{tr} 
{\bf Proof.}
In Section 3 
we have proved the Lemma \ref{key5} for any form ring with identity and 
shown that  the local-global principle is a consequence of Lemma \ref{key5}.
So, the result is true in particular if we have ${\E}(2n, R, \LMD)$ is a normal subgroup of ${\G}(2n, R, \LMD)$. 

To prove the converse we need $R$ to be finite as $k$-module, where $k$ is a 
commutative ring with identity ({\it i.e.}, a ring with trivial involution). 

Let $\alpha\in {\E}(2n,R,\LMD)$ and $\beta\in {\G}(2n,R, \LMD)$. Then 
$\alpha$ can be expressed as a product of matrices of the form 
$\vartheta_{ij}(\tn{ ring element})$ and 
$\vartheta_{i}(\tn{ column vector})$.
Hence we can write $\beta\alpha \beta^{-1}$ as a product of the matrices of the form 
$({\rm I}_{2n}+\beta\, \tn{M}(\star_1,\star_2)\beta^{-1})$, with $\langle \star_1, \star_2\rangle =0$, where 
$\star_1$ and $\star_2$ are suitably chosen standard basis vectors.
Now let $v=\beta \star_1$. Then we can write $\beta\alpha \beta^{-1}$ as a product of the matrices of the form 
$({\rm I}_{2n}+\beta\, \tn{M}(v,w)\beta^{-1})$, with $\langle v, w\rangle =0$ for some row vector $w$ in $R^{2n}$. 
We show that each $({\rm I}_{2n}+\tn{M}(v,w))\in {\E}(2n,R,\LMD)$.

Let $\gamma(X)={\rm I}_{2n}+X\tn{M}(v,w)$. Then $\gamma(0)={\rm I}_{2n}$. 
By Lemma \ref{HB} it follows that $S^{-1}R$ is a semilocal ring, where $S=k-\m$, 
$\m \in {\M}(k)$. Since $v\in {\rm Um}(2n, R)$, using Theorem \ref{swan} we get 
$$v\in {\E}(2n, S^{-1}R, S^{-1}\LMD) e_1,$$ 
hence $Xv\in {\E}(2n, S^{-1}R[X], S^{-1}\LMD[X]) e_1$. 
Therefore, applying Lemma \ref{key5} over 
$S^{-1}(A[X], \LMD[X])$ it follows that 
$$\gamma_{\m}(X)\in {\E}(2n,S^{-1}R[X],S^{-1}\LMD[X]).$$
Now applying Theorem \ref{LG}, it follows that 
$\gamma(X)\in {\E}(2n,R[X],\LMD[X])$.
Finally, putting $X=1$ we get the result. \hb

\section{Nipotent property for ${\k}$ of Hermitian groups}

We devote this section to discuss the study of nilpotent property of unstable ${\k}$-groups. 
The literature in this direction can be found in the work of 
A. Bak, N. Vavilov and R. Hazrat.
Throughout this section we assume $R$ is a commutative ring with identity, {\it i.e.}, we are 
considering trivial involution and $n\ge r+3$. Following is the statement of the theorem. 
\iffalse
In \cite{Bak}, Bak proved that for an almost 
commutative ring the unstable ${\k}$-group is nilpotent-by-abelian. 
Later in \cite{HV} Vavilov and Hazrat have generalized his result for 
Chevalley groups. Bak's construction uses a localization-completion method. In \cite{BBR}, the author with A. Bak and Ravi A. Rao showed that the localization part suffices to get the result if we restrict our concern for finite krull dimension. 
In their very recent paper ({\it cf.}~\cite{BHV}) Bak, Vavilov and Hazrat proved the relative case 
for the unitary and Chevalley groups. But, in my best knowledge, so far there is no result for Hermitian groups. I observe that using the above local-global principle,
arguing as in \cite{BBR}, it follows that the unstable ${\k}$ of Hermitian group is nilpotent-by-abelian. 
We emitting the proof of Theorem 4.1 in \cite{BBR}. 
Arguing in the same manner one can......... 
\fi

\begin{tr} \label{nil}
The quotient 
group $\frac{{\SH}(2n, R, a_1,\ldots, a_r)}{{\EH}(2n, R, a_1,\ldots, a_r)}$ is 
nilpotent for $n\ge r+3$. The  class of nilpotency is at the most \tn{max}  $(1, d+3-n)$,  where $d=\dim \,(R)$.
\end{tr}

The proof follows by emitting the proof of Theorem 4.1 in \cite{BBR}. 

\begin{lm} \label{sol3a} Let $I$ be an ideal contained in the Jacobson 
radical $J(R)$ of $R$, and $\beta\in {\SH}(2n, R, \LMD)$, with $\beta\equiv {\rm I}_n$ modulo $I$. Then there 
exists $\theta\in {\EH}(2n, R, a_1,\!\!\ldots\!\!, a_r)$ such that $\beta\theta$= the diagonal matrix 
$[d_1,d_2,\dots,d_{2n}]$, where each $d_i$ is a unit in $R$ with $d_i\equiv 1$ 
modulo $I$, and $\theta$ a product of 
elementary generators with each congruent to identity modulo $I$.
\end{lm} 
{\bf Proof.} The diagonal elements are units. Let $\beta=(\beta_{ij})$, where 
$d_i=\beta_{ii}=1+s_{ii}$ with $s_{ii}\in I\subset J(R)$, for 
$i=1,\ldots, 2n$, and $\beta_{ij}\in I\subset J(R)$ for $i\ne j$. 
First we make all the $(2n,j)$-th, and $(i,2n)$-th  entries zero, for $i=2,\ldots,n$, $j=2,\ldots,n$. 
Then repeating the above process we can reduce the size of $\beta$. 
Since we are considering trivial involution, we take 
$$\alpha \!\!
=\!\!\underset{j=1}{\overset{n}\Pi}hl_{nj}(-\beta_{2nj}d_{j}^{-1})\!\!
\underset{n+1\le j\le n+r}{\underset{n+r+1\le i\le 2n-1}\Pi} \!\!\! hm_i(-\zeta_jd_j^{-1})\!\!\!\!\!
\underset{n+r+1\le j\le 2n-1}{\underset{r+1\le i\le n-1} \Pi} \!\!\! h\eps_{in}(\beta_{\rho(n)\rho(i)}d_{j}^{-1}),$$
where $j=i-r$ and $\zeta_j=(0,\ldots,0,\beta_{2n  j})$, and \\
$$\gamma=\underset{r+1\le i\le 2n-1}{\underset{r+1\le j\le 2n-1}\Pi h\eps_{nj}(a_{i-r}(\star)d_{2n}^{-1})}
hr_n(\eta),$$ where $a_t=0$ for $t>r$, and 
$\eta=(\beta_{1 2n}d_{2n}^{-1},\beta_{2 2n}d_{2n}^{-1},\ldots,\beta_{n 2n}d_{2n}^{-1})$.
Then the last column and last row of $\gamma\beta\alpha$ become 
$(0,\dots,0, d_{2n})^t$, where $d_{2n}$ is a unit in $R$ and $d_{2n}\equiv 1$ modulo 
${\I}$. Repeating the process we can modify $\beta$ to the required form. \hb

\begin{pr} $(${\it cf.} Lemma 7, \cite{P}$)$
\label{sol5} Let $(R, \LMD)$ be a commutative form ring, {\it i.e.}, with trivial involution and $s$ be a non-nilpotent element in $R$ and $a\in R$. 
Then for $l\ge 2$  
 $$\left[\vartheta_{ij}\left(\frac{a}{s} \right), {\SH}(2n, s^lR)\right]
\subset {\EH}(2n, R).$$
More generally, 
$\left[\eps, {\SH}(2n, s^lR)\right]\subset {\EH}(2n, R)$, for $l\gg 0$ and $\eps\in {\EH}(2n, R_s)$. 
\end{pr}

{\bf Proof of Theorem \ref{nil}:}
Recall \vp 

 Let $G$ be a group. Define $Z^0=H$, $Z^1=[G,G]$ and 
$Z^i=[G,Z^{i-1}]$. Then $G$ is said to be nilpotent if $Z^r=\{e\}$ for
some $r>0$, where $e$ denotes the identity element of $H$.

Since the map ${\EH}(2n, R, a_1,\ldots, a_r)\ra {\EH}(2n, R/I, \ol{a}_1,\ldots, \ol{a}_r)$
is surjective we may and do assume that $R$ is a reduced 
ring. Note that if $n\ge d+3$, then the group 
${\SH}(2n, R, a_1,\ldots, a_r)/{\EH}(2n, R, a_1,\ldots, a_r)={\KH}_1(R, a_1,\ldots, a_r)$, 
which is abelian and hence nilpotent. So we consider the case $n\le d+3$. 
Let us first fix a $n$. We prove the theorem by induction 
on $d=\dim R$. Let $$G={\SH}(2n, R, a_1,\ldots, a_r)/{\EH}(2n, R, a_1,\ldots, a_r).$$ Let $m=d+3-n$ and 
$\alpha=[\beta, \gamma]$, for some $\beta\in G$ and $\gamma\in Z^{m-1}$. 
Clearly, the result is true for $d=0$. 
Let $\widetilde{\beta}$ be the pre-image of $\beta$ under the map 
$${\SH}(2n, R, a_1,\ldots, a_r)\ra {\SH}(2n, R, a_1,\ldots, a_r)/{\EH}(2n, R, a_1,\ldots, a_r).$$
If $R$ is local then arguing as Lemma \ref{sol3a} it follows that 
%$\EH_{2n}=\SH_{2n}$, hence 
we can choose a non-zero-divisor $s$ in $R$ such that  
$\widetilde{\beta}_s\in {\EH}(2n, R_s, a_1,\ldots, a_r)$.

Consider $\ol{G}$,
where bar denote reduction modulo $s^l$, for some $l\gg 0$. By the  
induction hypothesis $\ol{\gamma}=\{1\}$ in $\ol{{\SH}(2n, R)}$,  
where bar denote the reduction modulo the subgroup ${\EH}(2n, R)$. 
Since ${\EH}(2n, R)$ is a normal subgroup
of ${\SH}(2n, R)$, for $n\ge r+3$, 
by modifying $\gamma$ we may assume that
$\widetilde{\gamma}\in {\SH}(2n, R,s^lR, a_1,\ldots, a_r)$, where $\widetilde{\gamma}$ is the pre-image 
of $\gamma$ in ${\SH}(2n, R,a_1,\ldots, a_r)$. Now by Proposition \ref{sol5} it follows that 
$[\widetilde{\beta},\widetilde{\gamma}]\in {\EH}(2n, R,a_1,\ldots, a_r)$. 
Hence $\alpha=\{1\}$ in $G$. \hb 

\begin{re} \tn{In (\cite{BRK}, Theorem 3.1) it has been proved that 
the question of normality of the elementary subgroup and the local-global principle are 
equivalent for the elementary subgroups of the linear, symplectic and 
orthogonal groups over an almost commutative ring with identity. 
There is a gap in the proof of the statement $(3)\Ra (2)$ of Theorem 3.1 in 
\cite{BRK} (for an almost commutative ring). 
The fact that over a non-commutative semilocal ring 
the elementary subgroups of the classical groups 
acts transitively on the set of unimodular and istropic 
({\it i.e.}, $\langle\,v,v\rangle=0$) vectors of length $n\ge 2$ in the linear
case, and $n=2r\ge 4$ in the non-linear cases has been used 
in the proof, but it is not mentioned anywhere in the article. 
This was pointed by Professor R.G. Swan and he provided us a proof for the above result. }
\end{re} 
 {\bf Acknowledgment:} My sincere thanks to Professors R.G. Swan for giving me his permission to
reproduce his proof of Theorem \ref{swan} 
(he gave a proof for the symplectic and orthogonal groups as noted above). 
%(The proof of H. Bass' Lemma \ref{B} given here is a simple proof of it that 
%he used to give in his K-theory courses).
I thank DAE Institutes in India and ISI 
Kolkata for allowing me to use their infrastructure facilities in times. 
I am very much grateful to Prof. A. Bak and Prof. Nikolai Vavilov for their kind efforts to correct the manuscript, and I thank 
University of Bielefeld, NBHM, and IISER Pune for their financial supports for my visits. 
I thank Professors T.Y. Lam, D.S. Nagaraj, Ravi Rao, B. Sury and Nikolai Vavilov for many useful suggestions and editorial inputs.  
I would like to give some credit to Mr. Gaurab Tripathi for correcting few mathematical misprints.

{\small
 
\addcontentsline{toc}{chapter}{Bibliography}} 


\begin{thebibliography}{99} 
\bibitem{A} {\scshape E. Abe}; Chevalley groups over local rings, 
Tohoku Math. J. (2) 21 (1969), 474 -- 494. 
\bibitem{Bak1} {\scshape A. Bak}; ${\K}$-Theory of forms. Annals of Mathematics Studies, 
98. Princeton University Press, Princeton, N.J. University of Tokyo Press, 
Tokyo, (1981). 
\bibitem{Bak} {\scshape  A. Bak}; Nonabelian ${\K}$-theory: the nilpotent class of 
${\k}$ and general stability. ${\K}$-Theory  4  (1991),  no. 4, 363--397.  
\bibitem{Bak3} {\scshape A. Bak, G. Tang}; Stability for Hermitian ${\k}$. 
Journal of Pure and Applied Algebra 150 (2000), 107--121. 
\bibitem{Bak2} {\scshape A. Bak, V. Petrov, G. Tang}; Stability for Quadratic ${\k}$. 
${\K}$-Theory 29 (2003), 1--11. 
\bibitem{BBR} {\scshape A. Bak, R. Basu, R.A. Rao}; 
Local-global principle for transvection groups. Proceedings of The American Mathematical Society 138 (2010), no. 4, 
1191--1204.
\bibitem{BHV} {\scshape A. Bak, R. Hazrat, N. Vavilov}; Localization-completion strikes again: Relative ${\k}$ is 
nilpotent-by-abelian. Jounal of Pure and Applied Algebra 213 (2009), 1075--1085.
\bibitem{BV} {\scshape A. Bak, N. Vavilov}; Structure of hyperbolic unitary groups I, elementary subgroups. 
Alg. Colloquium 7 (2) (2000), 159--196.
\bibitem{Ba} {\scshape H. Bass}; Algebraic K-Theory, Benjamin, New York-Amsterdam, (1968). 
\bibitem{B} {\scshape H. Bass}; Unitary algebraic ${\K}$-theory.  Algebraic K-theory, III: 
Hermitian ${\K}$-theory and geometric applications (Proc. Conf., Battelle 
Memorial Inst., Seattle, Wash., 1972). Lecture Notes in Mathematics, 
Vol. 343, Springer, Berlin (1973), 57--265.
\bibitem{Bass} {\scshape H. Bass}; Quadratic modules over polynomial rings. Contribution to Algebra 
(Collection of papers dedicated to Ellis Kolchin) Academic Press, N.Y. (1977), 1--23. 
\bibitem{BRK} {\scshape R. Basu, R.A. Rao, R. Khanna}; 
On Quillen's local-global principle. Commutative Algebra and Algebraic 
Geometry (Bangalore, India, 2003), Contemp. Math. 390, American Mathematical Society, 
Providence, RI, (2005), 17--30.
\bibitem{BR} {\scshape R. Basu, R.A. Rao};  
Injective stability for ${\k}$ of the orthogonal group.  Journal of Algebra  323  (2010),  no. 2, 393--396.
\bibitem{BK} {\scshape A.J. Berrick, M.E. Keating}; An Introduction to Rings and Modules
with K-theory in view; Cambridge University Press, 65 (2000).
\bibitem{CR} {\scshape P. Chattopadhyay, R.A. Rao}; Elementary symplectic orbits and improved ${\k}$-stability. 
J. K-Theory 7,  no. 2 (2011), 389--403.
\bibitem{FRS} {\scshape J. Fasel, R.A. Rao, R.G. Swan}; On stably free modules over affine algebras. Publ. Math. Inst. Hautes \u{E}tudes Sci. 116 (2012), 223–-243.
\bibitem{Fu} {\scshape Fu An Li}; The structure of orthogonal groups over arbitrary commutative rings, 
Chinese Ann. Math. Ser. B 10 (1989), 341 -- 350. 
\bibitem{HR} {\scshape R. Hazrat}; Dimension theory and nonstable ${\k}$ of quadratic modules.  ${\K}$-Theory  27  
(2002),  no. 4, 293--328. 
\bibitem{HV} {\scshape R. Hazrat, N. Vavilov}; ${\k}$ of Chevalley groups are nilpotent.
Journal of Pure and Applied Algebra  179  (2003), no. 1-2, 99--116.
\bibitem{HV1} {\scshape R. Hazrat, N. Vavilov}; Bak's work on ${\K}$-theory of rings. On the occasion of his 65th 
birthday. Journal of ${\K}$-Theory 4, no. 1 (2009), 1--65.
\bibitem{HVZ} {\scshape R. Hazrat, N. Vavilov, Z. Zhang}; 
Relative unitary commutator calculus, and applications. J. Algebra 343  (2011), 107--137.
\bibitem{HSVZ} {\scshape R. Hazrat, A. Stepanov, N. Vavilov, Z. Zhang};
The yoga of commutators. Zap. Nauchn. Sem. S.-Peterburg. Otdel. Mat. Inst. Steklov. 
(POMI) 387 (2011), Teoriya Predstavlenii, Dinamicheskie Sistemy, Kombinatornye Metody. 
XIX, 53--82, 189; translation in J. Math. Sci. (N. Y.) 179, no. 6 (2011), 662--678. 
\bibitem{HO} {\scshape A.J. Hahn, O.T. O'Meara}; The Classical groups and ${\K}$-theory. With a foreword by 
J. Dieudonné. Grundlehren der Mathematischen Wissenschaften [Fundamental Principles of Mathematical Sciences], 
291. Springer-Verlag, Berlin, (1989).
\bibitem{J} {\scshape N. Jacobson}; Lectures on Quadratic Jordan algebras, Tata Istitute of Fundamental Research, 
Bombay, (1969). 
\bibitem{K1} {\scshape I.S. Klein, A.V. Mikhalev}; The Orthogonal Steinberg group over a ring with involution. 
(Russian)  Algebra i Logika  9  (1970) 145--166. 
\bibitem{K2} {\scshape I.S. Klein, A.V. Mikhalev}; The Unitary Steinberg group over a ring with involution. 
(Russian)  Algebra i Logika  9  (1970) 510--519.
\bibitem{KOP}  {\scshape V.I. Kopeiko}; The stabilization of Symplectic groups 
over a polynomial ring. Math. USSR. Sbornik 34 (1978), 655--669.
\bibitem{M} {\scshape K. McCrimmon}; A general theory of Jordan rings. Proc. Nat. Acad. Sci. U.S.A.  56  (1966).  
1072--1079.
\bibitem{P1} {\scshape R. Parimala}; Failure of Quadratic analog of Serre's Conjecture. Bulletin of 
American Mathematical Society 82 (1976b), 962--964. 
\bibitem{P2} {\scshape R. Parimala}; Failure of Quadratic analog of Serre's Conjecture. 
American Journal of Mathematics 100 (1978), 913--924.
\bibitem{P} {\scshape V.A. Petrov}; Odd unitary groups. 
%(Russian) Zap. Nauchn. 
%Sem. S.-Peterburg. Otdel. Mat. Inst. Steklov. (POMI) 305 (2003), Vopr. Teor. Predst. Algebr. 
%i Grupp. 10, 195--225, 241; translation in 
J. Math. Sci. (N. Y.) 130, no. 3 (2005), 4752--4766. 
\bibitem{PS} {\scshape V.A. Petrov, A.K. Stavrova}; Elementary subgroups in isotropic reductive groups. 
(Russian) Algebra i Analiz 20 (2008), no. 4, 160--188; 
translation in St. Petersburg Math. J. 20, no. 4 (2009), 625--644. 
\bibitem{RV} {\scshape R.A. Rao, W. van der Kallen}; Improved stability for 
${\k}$ and $\tn{WMS}_ d$ of a non-singular affine algebra. ${\K}$-theory 
(Strasbourg, 1992).  Asterisque  no. 226 (1994), 11, 411--420.
\bibitem{BRS} {\scshape R.A. Rao, R. Basu, Selby Jose}; Injective Stability for ${\k}$ 
of the Orthogonal group. Journal of Algebra 323 (2010), 393--396. 
\bibitem{SS} {\scshape Sergei Sinchuk}; Injective stability for unitary ${\k}$, revisited. (To appear). 
\bibitem{SUS} {\scshape A.A. Suslin}; On the structure of special Linear group over
polynomial rings. Math. USSR. Izv. 11 (1977), 221--238.
\bibitem{ST} {\scshape M.R. Stein};
Stability theorems for ${\k}$, ${\K}_2$ and related functors modeled on Chevalley groups. Japan. J. Math. (N.S.) 4 
no. 1 (1978),  77--108. 
\bibitem{SUSK}  {\scshape A.A. Suslin, V.I. Kopeiko}; Quadratic modules and Orthogonal 
groups over polynomial rings. Nauchn. Sem., LOMI 71 (1978), 216--250.
\bibitem{SUSV}  {\scshape A.A. Suslin, L.N. Vaserstein}; Serre's problem on projective modules over polynomial rings, 
and algebraic ${\K}$-theory. Izv. Akad SSSR. Ser. Mat. Tom 40, no. 5 (1976), 937--1001. 
\bibitem{T} {\scshape Guoping Tang}; Hermitian groups and $K$-theory.  
${\K}$-Theory  13, no. 3  (1998), 209--267.
\bibitem{Ta} {\scshape G. Taddei}; Normalite des groups elementaries dans les groupes de Chevalley 
sur an anneau, Application of Algebraic ${\K}$-Theory to Algebraic Geometry and Number Theory. Part II 
(Boulder, Colo., 1983), Contemp. Math., vol. 55, Amer. Math. Soc., Providence, RI, (1986), 693 -- 710.  
\bibitem{Tu} {\scshape M.S. Tulenbaev}; Schur multiplier of a group of elementary matrices of finite order, 
Zapiski Nauchnykh Seminarov Leningradskogo Otdeleniya Matematicheskogo Instuta im V.A. Steklova Akad. Nauk  
SSSR, Vol. 86 (1979), 162--169.
\bibitem{V} {\scshape L.N. Vaserstein}; On the Stabilization of the general Linear group 
over a ring.  Mat. Sbornik (N.S.) 79 (121) 405--424 (Russian); English 
translated in Math. USSR-Sbornik. 8 (1969), 383--400. 
\bibitem{V2} {\scshape L.N. Vaserstein}; Stabilization of Unitary and Orthogonal Groups 
over a Ring with Involution. Mat. Sbornik, Tom 81 (123) (1970) no. 3, 307--326.
\bibitem{V3} {\scshape L.N. Vaserstein}; Stabilization for Classical groups over rings. 
(Russian)  Mat. Sb. (N.S.)  93 (135)  (1974), 268--295, 327. 
\bibitem{vas2} {\scshape L.N. Vaserstein}; On the normal subgroups of ${\rm GL}\sb{n}$ over a ring. 
Algebraic $K$-theory, Evanston (1980) (Proc. Conf., Northwestern Univercisy, 
Evanston, Ill., (1980)), pp. 456--465, Lecture Notes in Mathematics, 854, Springer, 
Berlin-New York, (1981). 
\bibitem{We} {\scshape Weibe Yu}; 
Stability for odd unitary ${\k}$ under the $\LMD$-stable range condition. 
Journal of Pure and Applied Algebra, 217 (2013), 886–-891.
\bibitem{W} {\scshape J. Wilson}; The normal and subnormal structure of gneral Linear
groups. Proc. Camp. Phil. Soc. 71 (1972), 163--177. 
\end{thebibliography}
\end{document}